\documentclass[authoryear,preprint,11pt]{elsarticle}

\usepackage{amssymb}
\usepackage{amsthm}

\makeatletter
\def\ps@pprintTitle{%
 \let\@oddhead\@empty
 \let\@evenhead\@empty
 \def\@oddfoot{}%
 \let\@evenfoot\@oddfoot}
\makeatother

\journal{Applied Soft Computing}

\usepackage{algpseudocode}
\usepackage{algorithm}

 \newtheorem{Theorem}{Theorem}
 
 \newdefinition{Proposition}{Proposition}
 \newtheorem{Corollary}{Corollary}
 \newdefinition{Example}{Example}
 \newdefinition{Definition}{Definition}
 \newdefinition{Remark}{Remark}
 \newproof{Proof}{Proof}
 \newproof{pot}{Proof of Theorem \ref{thm2}}

\usepackage{amsmath}

\begin{document}

\begin{frontmatter}

\title{
A Fast Hybrid Primal Heuristic for Multiband Robust Capacitated Network Design with Multiple Time Periods\tnoteref{LabelTitle}
}

\tnotetext[LabelTitle]
{
This is the authors' final version of the paper published in Applied Soft Computing 26, 497-507, 2015, DOI: 10.1016/j.asoc.2014.10.016 \; .
The final publication is available at Elsevier ScienceDirect via http://dx.doi.org/10.1016/j.asoc.2014.10.016.
\\
Please note that this paper is an extended, improved journal version of the paper ``\emph{A hybrid primal heuristic for Robust Multiperiod Network Design}'' presented at EvoComNet2014 and published in: EvoApplications 2014: Applications of Evolutionary Computation, Granada, Spain, April 23-25, 2014, LNCS 8602, pp. 15-26, Springer, 2014. DOI: 10.1007/978-3-662-45523-4\_2.
\\
The paper was awarded the \emph{EvoComNet Best Paper Award} at \emph{EvoStar 2014}.
\\
\
\\
{
}
\\
}

\author[label1,label2,label3]{Fabio D'Andreagiovanni\corref{cor1}}
\ead{d.andreagiovanni@zib.de}
\ead[url]{http://www.dis.uniroma1.it/~fdag/}
\author[label1]{Jonatan Krolikowski}
\ead{krolikowski@zib.de}
\author[label1]{Jonad Pulaj}
\ead{pulaj@zib.de}

\address[label1]{
    Department of Optimization, Zuse-Institut Berlin (ZIB),
    Takustr. 7, 14195 Berlin, Germany
}

\address[label2]{
    DFG Research Center MATHEON, Technical University Berlin,
    Stra{\ss}e des 17. Juni 135, 10623 Berlin, Germany
}

\address[label3]{
    Einstein Center for Mathematics Berlin (ECMath),
    Stra{\ss}e des 17. Juni 135, 10623 Berlin, Germany
}

\cortext[cor1]{Corresponding author}

\begin{abstract}
We investigate the Robust Multiperiod Network Design Problem, a generalization of the Capacitated Network Design Problem (CNDP) that, besides establishing flow routing and network capacity installation as in a canonical CNDP, also considers a planning horizon made up of multiple time periods and protection against fluctuations in traffic volumes. As a remedy against traffic volume uncertainty, we propose a Robust Optimization model based on Multiband Robustness \citep{BuDA12a}, a refinement of classical $\Gamma$-Robustness by \cite{BeSi04} that uses a system of multiple deviation bands.
Since the resulting optimization problem may prove very challenging even for instances of moderate size solved by a state-of-the-art optimization solver, we propose a hybrid primal heuristic that combines a randomized fixing strategy inspired by ant colony optimization
and an exact large neighbourhood search.
Computational experiments on a set of realistic instances from the \cite{SNDlib} show that our original heuristic can run fast and produce solutions of extremely high quality associated with low optimality gaps.
\end{abstract}


\end{frontmatter}

\section{Introduction}

In the last two decades, telecommunications have increasingly pervaded our everyday life
and the volume of traffic sent and exchanged over networks has astonishingly increased: major companies like Nokia Siemens Networks expect that the increase in the amount of traffic will strongly continue, reaching a volume
of more than 1000 exabyte per year in fixed networks by 2015 \citep{Th09}.
This dramatic growth that telecommunications have experienced has greatly compounded the challenge for network professionals, who are facing design problems of increasing complexity and difficulty. In order to cope with traffic growth, the professionals have to plan much in advance how the network will be expanded in topology and capacity to accommodate the new traffic. This is especially important in the case of fixed networks, which require (costly) digging operations for the installation of cables in areas with possibly high population density.

To make the design task even more complicated, the future behaviour of traffic over a network is not exactly known when the network is designed and thus the decision problem is also affected by tricky data uncertainty. Until recent times, data uncertainty has been generally neglected in real studies. However, as indicated by recent industrial cooperations between industry and academia, (e.g., \citealt{BaEtAl14,BeKoNo11,BlDAKa13,KoKuRa13}), professionals are not only becoming aware of the importance of adopting mathematical optimization to take better decisions, but are also understanding the necessity of considering data uncertainty, in order to avoid unpleasant surprises like infeasibility of implemented solutions due to data deviations.

The task of designing a telecommunication network essentially consists in establishing the topology of the network and the technological features (e.g., transmission capacity and rate) of its elements, namely nodes and links.
One of the most studied problem in network design is the \emph{Capacitated Network Design Problem (CNDP)}: the CNDP consists in minimizing the total installation cost of capacity modules in a network to route traffic flows generated by users. The CNDP is a central problem in network optimization, which appears in many real-world applications. For an exhaustive introduction to it, we refer the reader to  \cite{AhMaOr93,Be98,BiEtAl98}.

In this paper,
we focus on the development of a new Robust Optimization model to tackle traffic uncertainty in the \emph{Multiperiod Capacitated Network Design Problem (MP-CNDP)}. This problem constitutes a natural extension of the classical CNDP, where, instead of a single design period, we consider
the design over a time horizon made up of multiple periods. Moreover, traffic uncertainty is taken into account to protect design solutions against deviations in the traffic input data, which may compromise feasibility and optimality of solutions.

We immediately stress that, though the problem of optimally designing networks over multiple time periods is not new and can be traced back at least to the seminal work by \cite{ChBr74}, to the best of our knowledge, the MP-CNDP has received very little attention and just a few works have investigated it - essentially, \cite{LaNaGe07} and \cite{GaRa12}.
Checking literature, several other works dealing with multiperiod design of networks can be found (just to make a couple of examples, \citealt{GeEtAl06} and \citealt{GuGr12}): however, all these works consider problems that are application-specific or are sensibly different from the more general setting that we consider here and we thus avoid a more detailed discussion of them (\cite{GeEtAl06} study capacity expansion problem in access networks with tree topology, whereas \cite{GuGr12} consider the design of utility networks modeled by non-linear mathematical programs).

Our main references for this work, namely \cite{LaNaGe07} and \cite{GaRa12},
point out the difficulty of solving multiperiod CNDP problems already for just two periods and even in (easier) contexts: \cite{LaNaGe07} consider the CNDP when traffic flows may be split, whereas \cite{GaRa12} consider a pure routing problem in satellite communications.
Our direct and more recent computational experience have confirmed the challenging nature of the MP-CNDP, even for instances of moderate size with a low number of time periods and solved by a state-of-the-art optimization solver.

Uncertain versions of the MP-CNDP where traffic uncertainty is considered have also been neglected, even though especially in the last years there has been an increasing interest in network design under traffic uncertainty for the single design period case (e.g., \citealt{BaEtAl14,BeKoNo11,KoKuRa13}).

\medskip
\noindent
In this work, our main original contributions are:
\begin{enumerate}
  \item the first Robust Optimization model for tackling traffic uncertainty in Multiperiod CNDP. Specifically, we adopt Multiband Robustness, a new model for Robust Optimization recently introduced by \cite{BuDA12a};
  \item a hybrid primal heuristic, based on the combination of a randomized rounding heuristic resembling \emph{ant colony optimization} \citep{DoDCGa99} with an exact large
    neighborhood search called RINS \citep{DaRoLP05}. We stress here that our aim was not to use a standard implementation of an ant colony algorithm: we wanted instead to strengthen the performance of the ant algorithm using highly valuable information from linear relaxations of the considered optimization problems. Using this information allowed us to define a very strong ant construction phase, which produces very high quality solutions already before the execution of any local search;
  \item analytically proving how to solve the linear relaxation of the Multiperiod CNDP in closed form, thus obtaining a substantial reduction in solution times w.r.t. our first algorithm presented in \citep{DAKrPu14};
  \item computational experiments over a set of realistic instances derived from the Survivable Network Design Library \citep{SNDlib}, showing that our hybrid algorithm is able to produce solutions of extremely high quality associated with very small optimality gap.
\end{enumerate}

\noindent
The remainder of this paper is organized as follows: in Section \ref{sec:CNDP}, we review a canonical model for the CNDP; in Section \ref{sec:MP-CNDP}, we introduce the Multiperiod CNDP and we study its linear relaxation;
in Section \ref{sec:Robust}, we introduce the new formulation for Robust Multiperiod CNDP; in Sections \ref{sec:ACO} and \ref{sec:computations}, we present our hybrid heuristic and computational results.

\section{The Capacitated Network Design Problem}
\label{sec:CNDP}

The CNDP is a central and highly studied problem in Network Optimization that appears in a wide variety of real-world applications (see \cite{AhMaOr93} and \cite{Be98} for an exhaustive introduction to it) and can be essentially described as follows:
given a network and a set of demands whose flows must be routed between vertices of the network, we want to install capacities on network edges and route the flows through the network, so that the capacity constraint of each edge is respected and the total cost of installing capacity is minimized.
More formally, we can characterize the CNDP through the following definition.
\begin{Definition}[The Capacitated Network Design Problem - CNDP]
Given
\begin{itemize}
  \item a network represented by a graph $G(V,E)$, where $V$ is the set of vertices and $E$ the set of edges,
  \item a set of commodities $C$, each associated with a traffic flow $d_{c}$ to route from an origin $s_{c}$ to a destination $t_{c}$,
  \item a set of admissible paths $P_c$ for routing the flow of each commodity $c$ from $s_{c}$ to $t_{c}$,
  \item a cost $\gamma_e$ for installing one module of capacity $\phi > 0$ on edge $e \in E$,
\end{itemize}
the CNDP consists in establishing the number of capacity modules installed on each edge $e \in E$ such that the resulting capacity installation has minimum cost and supports a feasible routing of the commodities. A feasible routing assigns each commodity $c \in C$ to exactly one feasible path $p \in P_c$.
\qed
\end{Definition}

\noindent
Referring to the notation introduced above and introducing the following two families of decision variables:
\begin{itemize}
  \item binary \emph{path assignment variables} $x_{cp} \in \{0,1\} \hspace{0.2cm} \forall \hspace{0.1cm} c \in C, p \in P_c$ such that:
    \begin{eqnarray*}
        x_{cp} &=& \left\{
                            \begin{array}{lll}
                                1 \hspace{0.2cm} \mbox{ if the entire traffic of commodity $c$ is routed on path $p$} \\
                                0 \hspace{0.2cm} \mbox{ otherwise,}
                            \end{array}
                        \right.
    \end{eqnarray*}
  \item integer \emph{capacity variables} $y_e \in \mathbb{Z}_{+}$, $\forall \hspace{0.1cm} e \in E$, representing the number of capacity modules installed on edge $e$,
\end{itemize}

\noindent
we can model the CNDP as the following \emph{integer linear program}:
\begin{align}
\min & \sum_{e \in E} \gamma_e \hspace{0.1cm} y_e
&&\mbox{(CNDP-IP)}
\nonumber
\\
&
\sum_{c \in C} \sum_{p \in P_c: \hspace{0.05cm} e \in p} d_{c} \hspace{0.1cm} x_{cp} \leq \phi \hspace{0.1cm} y_e
&&
e \in E
\label{CNDP_cnstr_capacity}
\\
&
\sum_{p \in P_c} x_{cp} = 1
&&
c \in C
\label{CNDP_cnstr_singlePath}
\\
&x_{cp} \in \{0,1\}
&& c \in C, p \in P_c
\nonumber
\\
&y_e \in \mathbb{Z}_{+}
&& e \in E \; .
\nonumber
\end{align}

\noindent
The objective function minimizes the total cost of capacity installation. Capacity constraints \eqref{CNDP_cnstr_capacity} impose that the summation of all flows routed through an edge $e \in E$ must not exceed the capacity installed on $e$ (equal to the number of installed modules represented by $y_e$ multiplied by the capacity $\phi$ granted by a single module). Constraints \eqref{CNDP_cnstr_singlePath} impose that the flow of each commodity $c \in C$ must be routed through a single path.

\begin{Remark}
This is an unsplittable version of the CNDP, namely the traffic flow of a commodity $c \in C$ \emph{cannot be split} over multiple paths going from $s_{c}$ to $t_{c}$, but must be routed on exactly one path.
Moreover, the set of feasible paths $P_c$  of each commodity is preset and constitutes an input of the problem.
This is in line with other works based on industrial cooperations (e.g., \citealt{BlGrWe07}) and with our experience (\citealt{BaEtAl14}), in which a network operator typically considers just a few paths that meet its own specific business and quality-of-service considerations and uses state-of-the-art routing schemes based on an Open Shortest Path First protocol.
\end{Remark}

\section{The Multiperiod Capacitated Network Design Problem}
\label{sec:MP-CNDP}

We define the multiperiod generalization of the CNDP by introducing a time horizon made up of a set of elementary periods $T = \{1,2, \ldots, |T|\}$. From a modeling point of view, this generalization requires to add a new index $t \in T$ to the decision variables, to represent routing and capacity installation decisions taken in each period. This is a simple modeling operation, which, however, greatly increases the size and complexity of the problem as pointed out in our computational section and by \cite{GaRa12} and \cite{LaNaGe07}.
The introduction of the new index leads to the new integer linear program:
\begin{align}
\min & \sum_{e \in E} \sum_{t \in T} \gamma_{e}^{t} \hspace{0.1cm} y_{e}^{t}
&&
\mbox{(MP-CNDP-IP)}
\nonumber
\\
&
\sum_{c \in C} \sum_{p \in P_c: \hspace{0.05cm} e \in p} d_{c}^{t} \hspace{0.1cm} x_{cp}^{t}
\leq \phi \hspace{0.1cm} \sum_{\tau = 1}^{t} y_{e}^{\tau}
&&
e \in E, t \in T
\label{MP-CNDP_robustcnstr_capacity}
\\
&
\sum_{p \in P_c} x_{cp}^{t} = 1
&&
c \in C, t \in T
\nonumber
\\
&x_{cp}^{t} \in \{0,1\}
&& c \in C, p \in P_c, t \in T
\nonumber
\\
&y_{e}^{t} \in \mathbb{Z}_{+}
&& e \in E, t \in T \; .
\nonumber
\end{align}

\noindent
Besides including the new index $t \in T$, the program presents a modified right-hand-side in the capacity constraint \eqref{MP-CNDP_robustcnstr_capacity}: for each time period $t$, we must consider the presence of all the capacity modules that are installed on an edge $e \in E$ from period $\tau = 1$ to period $\tau = t$.

Concerning costs and traffic demands, in what follows, we realistically assume that the cost per unit of capacity is non-increasing over time and that the demand associated with each commodity is non-decreasing over time (i.e., $\gamma_{e}^{t} \geq \gamma_{e}^{t+1}$, $\forall \hspace{0.1cm} t = 1,\ldots,|T|-1$ and $d_{c}^{t} \leq d_{c}^{t+1}$, $\forall \hspace{0.1cm} t = 1,\ldots,|T|-1$).

\subsection{\textbf{Computing the linear relaxation of MP-CNDP-IP}}

After having introduced the problem MP-CNDP-IP, we proceed to formally characterize in closed-form the value and optimal solution of its linear relaxation, namely the problem that we obtain if we relax the integrality requirements on variables $x$ and $y$ (so we have $0 \leq x_{cp}^{t} \leq 1$, $\forall c \in C, p \in P_c, t \in T$ and $y_{e}^{t} \geq 0$, $\forall e \in E, t \in T$). We refer to such relaxation as MP-CNDP-LP.

Characterizing the optimal value and the structure of an optimal solution of MP-CNDP-LP has been a crucial objective for us.
In our preliminary computational experience presented in the proceeding version of this paper \citep{DAKrPu14}, solving MP-CNDP-LP by using the state-of-the-art optimization solver IBM ILOG \cite{CPLEX} proved indeed slow. Using a solver like CPLEX to solve linear relaxations thus constituted a bottleneck that limited the possibility of running the heuristic a large number of times within the time limit.
By using instead the closed-form expression coming from our new theoretical results presented in this section, we were able to efficiently solve MP-CNDP-LP without using CPLEX, run the algorithm an incredibly higher number of times in the same time limit and gain a determinant speed up in the overall execution of the algorithm.

\bigskip

\noindent
In the case of an MP-CNDP including a single time period (i.e., $|T|=1$), it is well known that an optimal solution of the linear relaxation MP-CNDP-LP can be obtained by just considering the shortest path $p^{*} \in P_c$ for each commodity $c \in C$ and then by installing on each edge of $p^{*}$ the smallest number of capacity modules needed to support the traffic request $d_{c}$ of $c$ (see \citealt{AhMaOr93}).

Through the following propositions, we investigate how the optimal value and an optimal solution of the MP-CNDP-LP with $|T|>1$ look like, proving that they can be characterized efficiently. Note that in order to keep the exposition light, we decided to move the more complex proofs of the statements to the \ref{Appendix}.

\bigskip

\noindent
As first step, we characterize the relation between flows in consecutive time periods.

\begin{Proposition}
\label{Pr:consecutivePeriodFlow}
Let $\bar{x}^t_{c p}:=0$, $d^t_{c p}:=0$ and $\gamma^t_{e}:=\infty$ for $t=0$. There exists an optimal solution $\bar{x}$ of MP-CNDP-LP such that:
\begin{equation}
\label{eq:trafficRelation}
 d^{t-1}_c \hspace{0.1cm} \bar{x}^{t-1}_{c p}
\hspace{0.15cm} \leq \hspace{0.15cm}
d^t_c \hspace{0.1cm}  \bar{x}^t_{c p}
\hspace{1.0cm}
\forall \hspace{0.1cm} c \in C, p \in P_c, t \in T
\end{equation}
\end{Proposition}

\begin{Proof}
See the \ref{Appendix}.
\qed
\end{Proof}

\smallskip

\noindent
Using the previous Proposition, we are able to derive the following Theorem, which characterizes an optimal solution for the linear relaxation of MP-CNDP-IP. In particular, we can show that to get an optimal solution we just need to identify for each commodity $c \in C$ and time period $t \in T$ a shortest path in edge length $p_c^{*} \in P_c$ among the feasible paths $P_c$ of $c$. Of course, this operation can be done very fast and efficiently.

\begin{Remark}
We stress that in what follows, we will just focus on the determination of the optimal flows, since once we have an optimal flow, the capacity installation can be immediately derived over the network.
\end{Remark}

\begin{Theorem}
\label{Th:optSolMP-CNDP-LP}
Consider problem MP-CNDP-LP, namely the linear relaxation of MP-CNDP-IP, and let $p_c^{*} \in P_c$ be a shortest path in edge length for each commodity $c \in C$. An optimal solution of MP-CNDP-LP can be defined by routing in each time period $t \in T$ the entire flow of each commodity $c \in C$ on the shortest path $p_c^{*}$ and by installing on each edge of $p_c^{*}$ the exact capacity needed to route the traffic flow $d_{c}^{t}$.
\end{Theorem}

\begin{Proof}
See the \ref{Appendix}.
\qed
\end{Proof}

\smallskip

\noindent
The previous theorem efficiently characterizes an optimal solution for MP-CNDP-LP. However, our hybrid algorithm requires a feasible solution of the linear relaxation when the routing and capacity installation has been established for a number of consecutive time periods (see Section \ref{sec:ACO} for details). The following Corollary shows how to characterize a feasible solution for the modified linear relaxation and identifies the conditions under which this solution becomes optimal.

\begin{Corollary}
\label{Corol:optSolMP-CNDP-LP}
Consider problem MP-CNDP-IP and suppose that for time periods $t = 1, \ldots, \tau - 1$, $\tau \leq |T|$ and for all the commodities $c \in C$, feasible path assignments and capacity installations have been established (thus fixing the values of the corresponding variables $x$ and $y$). Moreover, suppose that for period $t = \tau$ feasible path assignments and capacity installations have been established for a subset of commodities $C' \subseteq C$.

Denote now by MP-CNDP-IP$^{FIX}$ the version of MP-CNDP-IP obtained from this variable fixing
and consider the corresponding linear relaxation MP-CNDP-LP$^{FIX}$.
Let $p_c^{*} \in P_c$ be a shortest path in edge cost for each commodity $c \in C$. A feasible solution for MP-CNDP-LP$^{FIX}$ can be defined by:
 \begin{itemize}
   \item routing in time period $t = \tau$ the entire flow of each $c \in C \setminus C'$ on the shortest path $p_c^{*}$ and by installing on each edge of $p_c^{*}$ the minimum number of capacity modules needed to route the traffic flow $d_{c}^{\tau}$;
\item routing in each time period $t = \tau +1, \ldots, |T|$ the entire flow of each $c \in C$ on the shortest path $p_c^{*}$ and by installing on each edge of $p_c^{*}$ the minimum number of capacity modules needed to route the traffic flow $d_{c}^{t}$;
 \end{itemize}
Furthermore, the solution determined above is optimal when for all the commodities $c \in C$, the path $p_c^{*}$ is chosen in every $t \in T$.
\end{Corollary}

\begin{Proof}
The feasibility of the solution built as specified above is clear. The optimality condition is instead a straightforward consequence of Theorem \ref{Th:optSolMP-CNDP-LP} and we omit the proof.
\qed
\end{Proof}

\noindent
The new theoretical results that we have introduced provide an alternative way to compute an (optimal) feasible solution of MP-CNDP-LP, which prove dramatically faster than the direct use of CPLEX. This allows us to greatly increase the number of executions of the ant construction phase in the same time limit.

\section{Robust Optimization for traffic-uncertain Multiperiod Network Design}
\label{sec:Robust}

After having introduced the multiperiod generalization MP-CNDP, we can proceed to consider its version that takes into account traffic uncertainty. To this end, in this section, we first state what we mean by traffic uncertainty, then we present fundaments of Robust Optimization, the methodology that we adopt to tackle data uncertainty, and finally we present a Robust Optimization model for the traffic-uncertain version of the MP-CNDP.

\subsection{\textbf{Traffic uncertainty and Robust Optimization}}

Uncertainty of traffic is naturally present in telecommunications network design, since the future behaviour of customers is not known in advance: the number of users and the traffic generated by them can just be estimated and these estimates can deeply differ from actual traffic conditions that will occur in the future (see \citealt{BaEtAl14}).
In what follows, we thus assume that the demands $d_c$ are \emph{uncertain} for all the commodities $c \in C$, i.e. their value is not known exactly when the optimization problem is solved.
In order to clarify the concept of traffic uncertainty, we anticipate here that we will model data uncertainty by a refined \emph{interval deviation model}. In an interval model, we assume to know a nominal value of traffic $\bar{d}_{c}$ and maximum negative and positive deviations $\delta_{c}^{-} \leq 0$, $\delta_{c}^{+} \geq 0$ from it. The (unknown) actual value $d_c$ is thus assumed to belong to the interval:
$$
\hspace{0.1cm}
d_c
\hspace{0.1cm} \in \hspace{0.1cm}
[\bar{d}_{c} + \delta_{c}^{-}, \hspace{0.2cm} \bar{d}_{c} + \delta_{c}^{+}] \;.
$$
In our direct experience with real-world network design, we have observed that professionals often identify $\bar{d}_{c}$ with the value of forecast traffic volume (e.g., an expected value derived from historical data), whereas the deviations $\delta_{c}^{-}, \delta_{c}^{+}$ are identified as the maximum deviations from the forecast considered relevant by the network designer, again using historical data as reference.

\begin{Example}[\textbf{traffic uncertainty}]
Consider two commodities $c_1, c_2$ associated with nominal traffic demands $\bar{d}_{c_1} = 200$ Mb, $\bar{d}_{c_2} = 300$ Mb and suppose that these values may deviate up to 10\%. So the maximum negative and positive deviations for $c_1, c_2$ are $\delta_{c_1}^{-} = -20$, $\delta_{c_1}^{+} = 20$ Mb, $\delta_{c_2}^{-} = -30$, $\delta_{c_2}^{+} = 30$ Mb, respectively. The actual values of traffic are therefore $d_{c_1} \in [180,220]$ Mb, $d_{c_2} \in [270,330]$ Mb.
\end{Example}

\noindent
As it is well-known from sensitivity analysis, dealing with data uncertainty in optimization problems may result really tricky: small variations in the value of input data may completely compromise the optimality and feasibility of produced solutions. Solutions that are supposed to be optimal may reveal to be heavily suboptimal, whereas solutions supposed to be feasible may reveal to be infeasible and thus meaningless when implemented. For a detailed discussion of the issues associated with data uncertainty, we refer the reader to \citep{BeElNe09,BeBrCa11}. The following example can immediately help to visualize the possibly catastrophic effects of neglecting data uncertainty.
\begin{Example}[\textbf{infeasibility caused by deviations}]
Consider again the commodities of Example 1 and suppose that in some link we have installed exactly the capacity to handle the sum of their nominal values (i.e., we have installed $200+300$ Mb of capacity).
This capacity installation neglects the fact that the demands may deviate up to 10\%. So, it is sufficient that just one commodity experiences a positive deviation to violate the capacity constraint of the link, thus making the design solution \emph{infeasible} in practice.
\end{Example}

\smallskip
\noindent
The previous example makes clear that we cannot afford to neglect traffic uncertainty and therefore risk that our design solution will turn out to be infeasible or of bad quality when implemented. As a consequence, we have decided to tackle data uncertainty by adopting a Robust Optimization (RO) approach. RO is a methodology for dealing with data uncertainty that has received a lot of attention in recent times and has been preferred to traditional methodologies like Stochastic Programming, especially thanks to its accessibility and computational tractability. We refer the reader to \cite{BeElNe09} and \cite{BeBrCa11} for an exhaustive introduction to theory and applications of RO and for a discussion about its determinant advantages over Stochastic Programming.

\noindent
RO is founded on two main facts:
\begin{itemize}
  \item the decision maker must define an \emph{uncertainty set}, which reflects his risk aversion and identifies the deviations of coefficients against which protection must be guaranteed;
  \item protection against deviations included in the uncertainty set is guaranteed through hard constraints that cut off all those feasible solutions that may become infeasible for some deviations of the uncertainty set.
\end{itemize}

\noindent
In a more formal way, suppose that we are given a generic integer linear program:
\begin{eqnarray*}
v
\hspace{0.1cm} = \hspace{0.1cm}
\max
\hspace{0.1cm}
c' \hspace{0.05cm} x
\hspace{0.5cm}
\mbox{ subject to }
\hspace{0.1cm}
x \in {\cal F}
=
        \{
        A \hspace{0.05cm} x \leq b,
        \hspace{0.15cm}
        x \in \mathbb{Z}^{n}_+
        \}
\end{eqnarray*}

\noindent
whose coefficient matrix $A$ is uncertain, namely we do not know the exact value of its entries. However, we can identify a family $\cal A$ of coefficient matrices that constitute possible valorizations of the uncertain matrix $A$. The family $\cal A$ is the uncertainty set of the robust problem. A \emph{robust optimal solution}, i.e. a solution that is protected against data deviations, can then be found by considering the following \emph{robust counterpart} of the original problem:
\begin{align*}
v^{{\cal R}}
\hspace{0.1cm} = \hspace{0.1cm}
\max
\hspace{0.1cm}
c' \hspace{0.05cm} x
\hspace{0.5cm}
\mbox{ subject to }
\hspace{0.1cm}
x \in {\cal R}
=
        \{
        \bar{A} \hspace{0.1cm} x \leq b
        \hspace{0.2cm}
        \forall \bar{A} \in {\cal A},
        \hspace{0.15cm}
        x \in \mathbb{Z}^{n}_+
        \}
\end{align*}

\noindent
The feasible set ${\cal R}$ of the robust counterpart includes only the solutions of ${\cal F}$ that are feasible for all the coefficient matrices in the uncertainty set ${\cal A}$. Consequently, ${\cal R}$ is a subset of the feasible set of the original problem, i.e. ${\cal R} \subseteq {\cal F}$.
We stress that the definition of robust counterpart can be extended to any mixed-integer linear program, involving continuous and integer variables at the same time.
Additionally, we remark that the decision maker should include in ${\cal A}$ coefficient matrices that reflect his specific risk aversion.

Ensuring protection against data deviations according to an RO paradigm comes at a price: the so-called \emph{price of robustness (PoR)} \citep{BeSi04}. The PoR is the deterioration of the optimal value of the robust counterpart w.r.t. the optimal value of the original problem (we have $v^{{\cal R}} \leq v$) and is caused by the hard constraints imposing robustness, which restrict the feasible set ${\cal F}$ to the (\emph{in general smaller}) set of robust solutions ${\cal R}$. The PoR depends upon the features of the uncertainty set: uncertainty sets reflecting higher levels of risk aversion of the decision maker will include more unlikely and extreme deviations, leading to higher protection yet associated with higher PoR; uncertainty sets reflecting low risk aversion will instead tend to neglect unlikely deviations, thus guaranteeing lower protection yet associated with lower PoR.

\begin{Example}[\textbf{protection against deviations}]
Following Example 2, a simple way to grant protection would be to install sufficient capacity to deal with the peak deviations of each commodity. So we should install 220+330 Mb of capacity. We note that in practice it is unlikely that all coefficients will experience the worst deviation, so one of the aim of ``smart'' RO models is to define appropriate uncertainty sets that result not too conservative, while guaranteeing a satisfying level of protection (for example, we could assume that at most one of the two demands will deviate from its nominal value).
\end{Example}

\noindent
In the next paragraph, we provide a description of the model of uncertainty that we adopt.

\subsection{\textbf{A concise introduction to Multiband Robust Optimization}}
\label{sec:multiband}

In this work, we tackle uncertainty through \emph{Multiband Robust Optimization} (MB), a new robust optimization model based on cardinality-constrained uncertainty set that was proposed by \cite{BuDA12a} and then extended and applied in a series of successive works (e.g, \citealt{BuDA12a,BuDA13,BuDA14,BuDARa14,BaEtAl14}). MB represents a refinement and generalization of the well-known $\Gamma$-Robustness ($\Gamma$-Rob) by \cite{BeSi04} that we developed to satisfy practical needs of our industrial partners in real-world applications (see \citealt{BuDA13,BaEtAl14}).

We recall here the main results of MB, referring  to the following generic \emph{uncertain} Mixed-Integer Linear Program (MILP):
\begin{align*}
\max
&
\sum_{j\in J} c_j \hspace{0.05cm} x_j
&&
(MILP)
\\
& \sum_{j\in J} a_{ij} \hspace{0.05cm} x_j \leq b_i
&& i\in I = \{1, \ldots, m\}
\\
& x_j \geq 0
&& j \in J = \{1, \ldots, n\}
\\
& x_j \in \mathbb{Z}_{+}
&& j \in J_{\mathbb{Z}} \subseteq J
\; .
\end{align*}

\noindent
where w.l.o.g we assume that the uncertainty only affects the coefficients $a_{ij}$ (uncertainty affecting cost coefficients or the right-hand-sides can be easily reformulated as coefficient matrix uncertainty).

Moving from the hypothesis that the actual value $a_{ij}$ of the coefficients is unknown,
the multiband uncertainty model at the basis of MB assumes that:
\begin{enumerate}
  \item for each coefficient $a_{ij}$, the decision maker knows its nominal value $\bar{a}_{ij}$
    as well as the maximum negative and positive deviations $d_{ij}^{K^{-}},d_{ij}^{K^{+}}$
    from $\bar{a}_{ij}$ (so $a_{ij} \in [\bar{a}_{ij}+d_{ij}^{K^{-}},\bar{a}_{ij}+d_{ij}^{K^{+}}]$);
  \item the overall single deviation band $[d_{ij}^{K-},d_{ij}^{K+}]$ of each coefficient $a_{ij}$ is partitioned into $K$ bands, defined on the basis of $K$ deviation values:
        $$
        -\infty<
        {d_{ij}^{K^{-}}<\cdots<d_{ij}^{-1}
        \hspace{0.1cm}<\hspace{0.2cm}d_{ij}^{0}=0\hspace{0.2cm}<\hspace{0.1cm}
        d_{ij}^{1}<\cdots<d_{ij}^{K^{+}}}
        <+\infty ;
        $$
  \item through these deviation values, $K$ deviation bands are defined, namely:
    a set of positive deviation bands $k\in \{1,\ldots,K^{+}\}$ and  a set of negative deviation bands $k \in \{K^{-}+1,\ldots,-1,0\}$, such that a band $k\in \{K^{-}+1,\ldots,K^{+}\}$ corresponds to the range $(d_{ij}^{k-1},d_{ij}^{k}]$, and band $k = K^{-}$
    corresponds to the single value $d_{ij}^{K^{-}}$;
  \item for each constraint $i \in I$ and each band $k\in K$, a lower bound $l_{ik}$ and
        an upper bound $u_{ik}$ on the number of deviations that may fall in $k$ are defined, so $0 \leq l_{ik} \leq u_{ik} \leq n$;
\item the number of coefficients that take
        their nominal value is not limited, i.e. $u_{i0}=n$ for all $i \in I$;
\item $\sum_{k \in K} l_{ik} \leq n$ for all $i \in I$, so that there always exists a feasible realization of the coefficient matrix.
\end{enumerate}

\noindent
We call this typology of uncertainty set a \emph{multiband uncertainty set}.

MB thus generalizes the uncertainty definition of the $\Gamma$-Rob model by \cite{BeSi04}: the single deviation band is partitioned into multiple bands and each band $k \in K$ is associated not only with an upper bound $u_{ik}$, but also with a lower bound $l_{ik}$ on the number of coefficients deviating in that band. The lower bound improves the modeling power of the decision maker and, more importantly, allows to take into account the presence of negative value deviations, that are neglected in $\Gamma$-Rob. Of course, taking into account negative deviations allows to improve the modeling of deviations commonly found in real-world problems and, even more critically, reduces the value of the overall worst deviation and thus the price of robustness. The multiband model results particularly suitable to model histograms that are commonly adopted by practitioners to visualize and analyze data deviations (see \citealt{BuDA12a,BaEtAl14}).

Since the robust optimization paradigm entails that we must be protected against any possible deviation considered in the uncertainty set, the robust counterpart of MILP under multiband uncertainty is:
\begin{align}
\max
&
\sum_{j\in J} c_j \hspace{0.05cm} x_j
&&
\label{NonLinearRobustCounterpart}
\\
& \sum_{j\in J} \bar{a}_{ij} \hspace{0.05cm} x_j + DEV_i(x,MB) \leq b_i
&& i\in I
\nonumber
\\
& x_j \geq 0
&& j \in J
\nonumber
\\
& x_j \in \mathbb{Z}_{+}
&& j \in J_{\mathbb{Z}} \subseteq J
\nonumber
\; ,
\end{align}

\noindent
where an additional term $DEV_i(x,MB)$ is introduced in every feasibility constraint to represent the maximum total deviation that could be incurred by constraint $i$ under the multiband uncertainty set for a solution $x$. This problem is actually non-linear since the term $DEV_i(x,MB)$ hides the following maximization problem:
\begin{align}
DEV_i(x,MB) = \max & \sum_{j\in J} \sum_{k\in K} d_{ij}^{k} \hspace{0.1cm} x_j \hspace{0.1cm} y_{ij}^{k} &&
\label{MBdev_objFunction}
\\
& l_{ik} \leq \sum_{j\in J} y_{ij}^{k} \leq u_{ik}
&& k\in K
\label{MBdev_constraint}
\\
& \sum_{k\in K} y_{ij}^{k} \leq 1
&& j \in J
\label{MBgub_constraint}
\\
& y_{ij}^{k} \in \{0,1\}
&& j \in J, k\in K
\label{MBgub_variables}
\; .
\end{align}

\noindent
A binary variable $y_{ij}^{k}$ of the problem is equal to $1$ when the coefficient $j$ in constraint $i$ falls in deviation band $k$ and is equal to $0$ otherwise. Each coefficient $j$ in the constraint must fall in at most one deviation band, thus requiring the introduction of the family of constraints \eqref{MBdev_constraint} (note that when $\sum_{k\in K} y_{ij}^{k} = 0$ then it is like having the coefficient falling in the zero deviation band). Constraints \eqref{MBgub_constraint} impose the bounds on the number of coefficients that may deviate in each band $k \in K$. Finally, the objective function \eqref{MBdev_objFunction} aims at maximizing the deviation allowed by the multiband uncertainty set for a given solution $x$ and constraint $i$.

The robust counterpart \eqref{NonLinearRobustCounterpart} is a (non-linear) max-max problem, since it includes the binary program (\ref{MBdev_objFunction}-\ref{MBgub_variables}). This is anyway not a real issue, since we have proved that the robust counterpart is equivalent to a \emph{compact and linear} mixed-integer linear program, as stated in the following theorem.
\begin{Theorem}[\textbf{B\"using and D'Andreagiovanni 2012}]\label{th:compact_MILP}
The robust counterpart of problem (MILP) for the multiband uncertainty set $MB$ is equivalent to the following compact and linear mixed-integer program:
\begin{align}
\max & \sum_{j\in J} c_{j} \hspace{0.1cm} x_{j}
&&\mbox{(Rob-MILP)}
\label{robustCompactMILP}
\\
&\sum_{j\in J} \bar{a}_{ij} \hspace{0.1cm} x_{j}
+ \sum_{k\in K} \theta_{k} \hspace{0.1cm} w_{i}^{k}
+ \sum_{j\in J} z_{ij} \leq b_{i}
&&
i\in I
\label{robustMILP_modConstraint}
\\
&w_{i}^{k} + z_{ij} \geq d_{ij}^{k} \hspace{0.1cm} x_{j}
&&
i \in I, j\in J, k\in K
\label{robustMILP_dualConstr1}
\\
&w_{i}^{k}  \in \mathbb{R}
&&i \in I, k\in K
\label{robustMILP_dualVar1}
\\
&z_{ij} \geq 0
&&i \in I,j\in J
\label{robustMILP_dualVar2}
\\
&x_{j} \geq 0
&&j\in J
\nonumber
\\
&x_j \in \mathbb{Z}_{+}
&&j \in J_{\mathbb{Z}} \subseteq J \; .
\nonumber
\end{align}
\end{Theorem}

\noindent
This problem includes $K \cdot m + n \cdot m$ additional continuous variables (\ref{robustMILP_dualVar1}-\ref{robustMILP_dualVar2}) and $K \cdot n \cdot m$ additional constraints \eqref{robustMILP_dualConstr1} to clear up the non-linearity of the trivial robust counterpart \eqref{NonLinearRobustCounterpart}. Moreover, constraints \eqref{robustMILP_modConstraint} include additional terms and involve the values $\theta_{k} \geq 0$ that express the number of deviations that occur in each deviation band $k \in K$ (these values constitute the so-called \emph{profile} of the multiband uncertainty set and are derived on the basis of the bounds $l_k,u_k$ - see \cite{BuDA13} for details). The resulting formulation has thus the nice properties of being \emph{compact} and \emph{linear}.
The proof of the theorem is based on pointing out the integrality of the polyhedron associated with and on exploiting strong duality. We refer the reader to \cite{BuDA12a,BuDA13} for the formal complete statement and proof of the presented theorem.

Theorem \eqref{th:compact_MILP} is a central result from the theory of Multiband Robust Optimization that we will use to derive a robust model for the MP-CNDP.

\subsection{\textbf{Multiband-Robust Multiperiod Network Design}}

We proceed now to use Multiband Robust Optimization and the related Theorem \eqref{th:compact_MILP} to tackle traffic uncertainty in the MP-CNDP.
If we denote by ${\cal D}$ the uncertainty set associated with the demands of the commodities, we can write the general form of the robust counterpart of the MP-CNDP-IP as follows:
\begin{align}
\min & \sum_{e \in E} \sum_{t \in T} \gamma_{e}^{t} \hspace{0.1cm} y_{e}^{t}
&&
\label{NonLinearRobustCounterpart_MPCNDP}
\\
&
\sum_{c \in C} \sum_{p \in P_c: \hspace{0.05cm} e \in p} \bar{d}_{c}^{t} \hspace{0.1cm} x_{cp}^{t}
+ DEV_{e}^{t}(x, {\cal D})
\leq \phi \sum_{\tau = 1}^{t} y_{e}^{\tau}
&&
e \in E, t \in T
\label{Rob-CNDP_robustcnstr_capacity}
\\
&
\sum_{p \in P_c} x_{cp}^{t} = 1
&&
c \in C, t \in T
\nonumber
\\
&x_{cp}^{t} \in \{0,1\}
&& c \in C, p \in P_c, t \in T
\nonumber
\\
&y_{e}^{t} \in \mathbb{Z}_{+}
&& e \in E, t \in T \; .
\nonumber
\end{align}

\noindent
This robust counterpart differs from the MP-CNDP-IP in the capacity constraints \eqref{Rob-CNDP_robustcnstr_capacity}: these constraints indeed consider the \emph{nominal} traffic demands values $\bar{d}_{c}^{t}$ and include the terms $DEV_{e}^{t}(x, {\cal D})$ to represent the total maximum positive deviation that demands may experience on edge $e$ in period $t$ for a routing vector $x$ and the uncertainty set ${\cal D}$.

We structure the uncertainty set ${\cal D}$ according to the principles of Multiband Robust Optimization introduced in the previous subsection. Specifically, we build a \emph{multiband uncertainty set} for the MP-CNDP as follows:
\begin{enumerate}
  \item for each commodity $c \in C$ and time period $t \in T$, we know the nominal value $\bar{d}_{c}^{t}$ of the traffic coefficient and maximum negative and positive deviations $\delta_{c}^{t-} \leq 0$, $\delta_{c}^{t+} \geq 0$ from it. The actual value $d_{c}^{t}$ is then such that $d_{c}^{t} \in [\bar{d}_{c}^{t} + \delta_{c}^{t-}, \hspace{0.1cm} \bar{d}_{c}^{t} + \delta_{c}^{t+}]$;
  \item the overall deviation range $[\bar{d}_{c}^{t} + \delta_{c}^{t-}, \hspace{0.1cm} \bar{d}_{c}^{t} + \delta_{c}^{t+}]$ of each coefficient $d_{c}^{t}$ is partitioned into $K$ bands, defined on the basis of $K$ deviation values:
      \\
{\small        $
        -\infty<
        {\delta_{c}^{t-} = \delta_{c}^{tK^{-}}<\cdots<\delta_{c}^{t-1}
        < \delta_{c}^{t0}=0<
        \delta_{c}^{t1}<\cdots<\delta_{c}^{tK^{+}}} = \delta_{c}^{t+}
        <+\infty ;
        $
}
  \item through these deviation values, $K$ deviation bands are defined, namely:
    a set of positive deviation bands $k\in \{1,\ldots,K^{+}\}$ and  a set of negative deviation bands $k \in \{K^{-}+1,\ldots,-1,0\}$, such that a band $k\in \{K^{-}+1,\ldots,K^{+}\}$ corresponds to the range $(d_{c}^{tk-1},d_{c}^{tk}]$, and band $k = K^{-}$
    corresponds to the single value $d_{c}^{tK^{-}}$;
  \item for each capacity constraint $\eqref{Rob-CNDP_robustcnstr_capacity}$ defined for an edge $e \in E$ and period $t \in T$ and for each band $k\in K$, we introduce two values $l_{e}^{tk},u_{e}^{tk}:$ $0 \leq l_{e}^{tk} \leq u_{e}^{tk} \leq n_{e}^{t}$ to represent the lower and upper bound on the number of traffic coefficients whose value deviates in band $k$ ($n_{e}^{t}$ is the number of uncertain coefficients in the constraint). These bounds can be used to derive the profile of the uncertainty set, namely the values $\theta_{e}^{tk} \geq 0$ indicating the exact number of coefficients deviating in each band (see \cite{BuDA13} for details about the definition of profile);
  \item the number of coefficients that take
        their nominal value is not limited, i.e. $u_{e}^{t0}=n_{e}^{t}$ for all $e \in E,t \in T$;
  \item $\sum_{k \in K} l_{e}^{tk} \leq n_{e}^{t}$ for all $e \in E,t \in T$, so that there always exists a feasible realization of the coefficient matrix.
\end{enumerate}

\noindent
By using the previous characterization of the multiband uncertainty set and Theorem \eqref{th:compact_MILP}, we can reformulate the non-linear robust counterpart \eqref{NonLinearRobustCounterpart_MPCNDP} as the following \emph{linear and compact} robust counterpart:
\begin{align}
\min & \sum_{e \in E} \sum_{t \in T} \gamma_{e}^{t} \hspace{0.1cm} y_{e}^{t}
&&
\mbox{(Rob-MP-CNDP)}
\nonumber
\\
&
\sum_{c \in C} \sum_{p \in P_c: \hspace{0.05cm} e \in p} \bar{d}_{c}^{t} \hspace{0.1cm} x_{cp}^{t}
\hspace{0.1cm} +
&&
\nonumber
\\
&
+
\sum_{k \in K} \theta_{e}^{tk} \hspace{0.05cm} w_{e}^{tk}
+ \sum_{c \in C} \sum_{p \in P_c: \hspace{0.05cm} e \in p} z_{ecp}^{t}
\leq
\phi \sum_{\tau = 1}^{t} y_{e}^{\tau}
&&
e \in E, t \in T
\nonumber
\\
&
w_{e}^{tk} + z_{ecp}^{t} \geq \delta_{c}^{tk} \hspace{0.1cm} x_{cp}^{t}
&&
e \in E, c \in C, p \in P_c: e \in p,
\label{RobCNDP_cnstr_dual}
\\
&
&&
t \in T, k \in K
\nonumber
\\
&w_{e}^{tk} \in \mathbb{R}
&& e \in E, t \in T, k \in K
\label{RobCNDP_dual1}
\\
&z_{ecp}^{t} \geq 0
&& e \in E, c \in C, p \in P_c: e \in p, t \in T
\label{RobCNDP_dual2}
%
\\
&
\sum_{p \in P_c} x_{cp}^{t} = 1
&&
c \in C, t \in T
\nonumber
\\
&x_{cp}^{t} \in \{0,1\}
&& c \in C, p \in P_c, t \in T
\nonumber
\\
&y_{e}^{t} \in \mathbb{Z}_{+}
&& e \in E, t \in T \; .
\nonumber
\end{align}

\noindent
This formulation includes the additional constraints \eqref{RobCNDP_cnstr_dual} and variables \eqref{RobCNDP_dual1}, \eqref{RobCNDP_dual2} to linearly reformulate the original (non-linear) problem including the term $DEV_{e}^{t}(x, {\cal D})$ in each capacity constraint.
Rob-MP-CNDP is the problem that we want to solve in the computational section to get robust solutions to the MP-CNDP.

\section{A hybrid primal heuristic for the Rob-MP-CNDP}
\label{sec:ACO}

In principle, we can get a robust optimal solution for Rob-MP-CNDP by using any commercial mixed-integer programming software, such as CPLEX. However, as showed in the computational section, solving Rob-MP-CNDP constitutes a difficult task even when considering a small number of time periods and using a state-of-the-art solver like CPLEX: after several hours of computation, solutions are still typically of low quality and far away from the optimum.
As a remedy, we were attracted by the effectiveness of (hybrid) MIP-based and bio-inspired heuristics in hard network design problems. Valid examples of such effectiveness are provided by \cite{CrGe02}, proposing a cooperative parallel tabu search algorithm for the single-period CNDP tested on transshipment networks, by \cite{DeDAKa13}, proposing a linear relaxation-based decomposition method for a variant of the single-period CNDP related to fair routing in wireless mesh networks, and \cite{KlEtAl12}, proposing a multiobjective evolutionary algorithm to solve a single-period CNDP arising in the design of large communications networks.
Concerning hybrid heuristics, we refer the reader to \cite{BlEtAl11} for a recent survey.
In the case of Rob-MP-CNDP, we developed a fast hybrid primal heuristic based on the combination of a randomized variable-fixing strategy resembling Ant Colony Optimization (ACO) and an exact large neighbourhood search.

It is widely known that ACO is a metaheuristic that was inspired by the foraging behaviour of ants. The seminal work by \cite{DoMaCo96} presenting an ACO algorithm for combinatorial problems has later been extended to integer and continuous problems (e.g., \citealt{DoDCGa99}) and has been followed by hundreds of other papers proposing refinements of the basic algorithms (e.g., \citealt{GaMoWe12,Ma99}) and investigating applications to relevant optimization problems (see \citealt{Bl05} for an overview).
An ACO algorithm presents the general structure specified in Algorithm \ref{(Gen-ACO)}. A loop consisting of two phases is executed until an arrest condition is satisfied: in a first phase, an ant builds up a solution under the guidance of probabilistic functions of variable fixing that resemble pheromone trails; then the pheromone trails are updated on the basis of how effective the adopted variable fixing has resulted. Once that the arrest condition is reached, a daemon action phase takes place and some solution improvement strategy is applied to bring the feasible solution built by the ants to a (local) optimum.
\begin{algorithm}
\caption{General structure of an ACO algorithm (Gen-ACO)}
\label{(Gen-ACO)}
\begin{algorithmic}[1]
\While{an arrest condition is reached}
    \State ant-based solution construction
    \State pheromone trail update
\EndWhile
\State daemon actions
\end{algorithmic}
\end{algorithm}

We now proceed to detail each phase of the previous sketch for our hybrid ACO-exact algorithm for the Rob-MP-CNDP. Our algorithm is defined hybrid since after having passed the ACO construction phase, the daemon action phase operates an exact large neighborhood search that is formulated as an integer linear program and is solved \emph{exactly} relying on the power of modern state-of-the-art commercial solvers.

\subsection{Ant-based solution construction.}
In the first step of the inner cycle of (Gen-ACO), a number $m \geq 0$ of \emph{ants} are defined and each ant iteratively builds a feasible solution for the problem. At a generic iteration of the construction, the ant is in a \emph{state} corresponding with a \emph{partial solution} and can make a further step towards completing the solution by making a \emph{move}: a move corresponds to fixing the value of some not-yet-fixed variable. The move that is executed is chosen according to a probability function: the function specifies the probability of implementing a move and thus of fixing a variable and is derived by combining an \emph{a priori} measure of the efficacy of the move and an \emph{a posteriori} measure of the efficacy of the move. More in detail, the probability of choosing a move that fixes a variable $j$ after having fixed a variable $i$ is specified by the following canonical formula:
\begin{equation}
\label{canonicalProbACO}
p_{ij} = \frac{\tau_{ij}^{\beta} + \eta_{ij}^{\delta}}
                {\sum_{f \in F} \tau_{if}^{\beta} + \eta_{if}^{\delta}}
\; \;,
\end{equation}
where $\tau_{ij}$ represents the a priori measure of efficacy, commonly called \emph{pheromone trail} value, and $\eta_{ij}$ represents the a posteriori measure of efficacy, commonly called \emph{attractiveness}. Note that this formula is influenced by the two parameters $\beta, \delta$ that appear as exponents of the measures and should be chosen by the decision maker on the basis of the specific problem considered.
For a more detailed description of the elements and actions of this phase, we refer the reader to the paper \citep{Ma99}, which presents ANTS, an improved ant algorithm used for solving the quadratic assignment problem, which we have taken as reference for our work. We considered ANTS particularly attractive as it proposes a series of refinements for classical ACO that allow to better exploit polyhedral information about the problem. Specifically, \citep{Ma99} sketches some ideas about how alternative formulations of the original problem could be exploited to define the \emph{pheromone trail} and the \emph{attractiveness} values. As additional desirable feature, ANTS makes use of a reduced number of parameters and adopts more efficient mathematical operations w.r.t. the canonical ant algorithms (products instead of exponentiations). For an exhaustive description of ANTS, we refer the reader to the paper by \cite{Ma99}.

Before describing how our ANTS implementation is structured, we make some preliminary considerations. The formulation Rob-MP-CNDP is based on four families of variables: 1) the path assignment variables $x_{cp}^{t}$; 2) the capacity variables $y_{e}^{t}$; 3-4) the auxiliary variables $w_{e}^{tk}$, $z_{ecp}^{t}$ coming from robust dualization.
Though we have to deal with four families, we can notice that once routing decisions are taken over the entire time horizon, then we can immediately derive the capacity installation of minimum cost.
Indeed, once the values of all path assignment variables are fixed, the routing is completely established and the worst traffic deviation term $DEV_{e}^{t}(x, {\cal D})$ can be efficiently derived without need of using the auxiliary variables $w_{e}^{tk}$, $z_{ecp}^{t}$ (this can be indeed traced back to solving a min-cost flow problem, as explained in \citealt{BuDA12a,BuDA13}). So we can easily derive the total traffic $D_{e}^{t}$ sent over an edge $e$ in period $t$ in the worst case. We can then derive the minimum cost installation by a sequential evaluation from period 1 to period $|T|$, keeping in mind that we must have $\left\lceil \frac{D_{e}^{t}}{\phi} \right\rceil$ capacity modules on $e$ in $t$ to accommodate the traffic.
In the ant-construction phase, we can consequently limit our attention to the binary assignment variables and we introduce the concept of \emph{routing state}.
\begin{Definition}[Routing state - RS]
  Let $P = \bigcup_{c \in C} P_c$ and let $R \subseteq C \times P \times T$ be the subset of triples $(c,p,t)$ representing the assignment of path $p \in P_c$ to commodity $c \in C$ in period $t \in T$. A \emph{routing state} is an assignment of paths to a subset of commodities in a subset of time periods which excludes that multiple paths are assigned to a single commodity. Formally:
        $$
        RS \subseteq R:
       \hspace{0.2cm}
       \not \exists (c_1,p_1,t_1), (c_2,p_2,t_2) \in RS:
       \hspace{0.1cm}
       c_1 = c_2
       \hspace{0.1cm} \wedge \hspace{0.1cm}
       p_1,p_2 \in P_{c_1}
     \hspace{0.1cm} \wedge \hspace{0.1cm}
       t_1 = t_2 \; .
       $$
\end{Definition}

\noindent
We say that a routing state RS is \emph{complete} when it specifies the path used by each commodity in each time period (thus $|RS| = |C||T|$). Otherwise the RS is called \emph{partial} and we have $|RS| < |C||T|$).

In the ANTS algorithm that we propose, we decided to assign paths considering time periods and commodities in a pre-established order. Specifically, we establish the routing in each time period separately, starting from $t=1$ and continuing up to $t=|T|$, and in each time period commodities are sorted in descending order w.r.t. their nominal traffic demand. Formally, this is operated through the cycle specified in Algorithm \ref{ALGcompleteRouting} that builds a \emph{complete routing state}.
\begin{algorithm}
\caption{Construction of a complete routing state}
\label{ALGcompleteRouting}
\begin{algorithmic}[1]
\For{$t := 1$ to $|T|$}
    \State sort $c \in C$ in descending order of $\bar{d}_{ct}$
    \For{(sorted $c \in C$)}
    \State assign a single path $p \in P_c$ to $c$
    \EndFor
\EndFor
\end{algorithmic}
\end{algorithm}

\medskip
\noindent
For an iteration $(t,c)$ of the nested cycles of Algorithm \ref{ALGcompleteRouting}, the assignment of a path to a commodity is equivalent to an ant that moves from a partial routing state $RS_i$ to a partial routing state $RS_j$ such that:
$$
RS_j = RS_i \cup \{(c,p,t)\} \hspace{0.1cm} \mbox{ with } p \in P_c\;  .
$$
We note that, by the definition of routing state, a sequence of moves is actually a sequence of fixings of decision variables, as in \citep{Ma99}.

The probability that an ant $k$ moves from a routing state $i$ to a more complete routing state $j$, chosen among a set of feasible routing states, is defined by the improved formula proposed by \cite{Ma99}:
$$
p_{ij}^{k} = \frac{\alpha \hspace{0.1cm} \tau_{ij} + (1-\alpha) \hspace{0.1cm} \eta_{ij}}
                {\sum_{f \in F} \alpha \hspace{0.1cm} \tau_{if} + (1-\alpha) \hspace{0.1cm} \eta_{if}}
\; ,
$$
where $\alpha \in [0,1]$ is a parameter assessing the relative importance of trail and attractiveness. This formula presents two peculiar advantages over the canonical formula \eqref{canonicalProbACO}: it adopts the single parameter $\alpha$ in place of the two parameters $\beta, \delta$ and uses $\alpha$ as coefficient of a product instead of an index of an exponentiation.

As discussed in \cite{Ma99}, the trail values $\tau_{ij}$  and the attractiveness values $\eta_{ij}$ should be provided by suitable lower bounds of the considered optimization problem. In our particular case:
1) $\tau_{ij}$ is derived from the values of the variables in the solution associated with the linear relaxation of the robust counterpart Rob-MP-CNDP;
2) $\eta_{ij}$ is equal to the value of a (good) feasible solution of the linear relaxation MP-CNDP-LP where a subset of the decision variables have fixed values as effect of the fixing decision taken in previous steps of the algorithm.

\subsection{Daemon actions: Relaxation Induced Neighborhood Search}
At the end of the ant-construction phase, we attempt at improving the quality of the feasible solution found by executing an \emph{exact local search} in a \emph{large neighborhood}. In particular, we adopt a modified \emph{Relaxation Induced Neighborhood Search} (RINS) (see \cite{DaRoLP05} for an exhaustive description of the method).

In an integer linear program, we look for a solution that is integral and at the same time guarantees the best objective value. In RINS, it is observed that an optimal solution of the linear relaxation of the problem provides an objective value that is better than that of any feasible solution. However, at the same time, such an optimal solution is fractional and therefore does not guarantee integrality. On the contrary, a feasible solution guarantees integrality, but provides a worse objective value. The fact that a variable is fixed to the same value in the optimal solution of the linear relaxation and in a feasible solution is a good indication that the fixing of the variable is good and should thus be maintained. Given a feasible solution, RINS profits from these observations by defining search neighborhoods, where those variables that have the same value in the feasible solution and in the optimal solution of the linear relaxation are fixed, while the other are free to vary their value. The neighborhood is then explored exhaustively by formulating the search as an integer linear program, which is solved exactly, possibly including an arrest condition such as a solution time limit.

In our specific case, let $(\bar{x},\bar{y})$ be a feasible solution of Rob-MP-CNDP found by an ant and $(x^{\small LR},y^{\small LR})$ be an optimal (continuous) solution of the linear relaxation of Rob-MP-CNDP.
Our modified RINS \emph{(mod-RINS)} entails to solve exactly through an optimization solver like CPLEX a subproblem of Rob-MP-CNDP where:
\begin{enumerate}
    \item we fix the variables $x$ whose value in $(\bar{x},\bar{y})$ and $(x^{\small LR},y^{\small LR})$ differs of at most $\epsilon > 0$, i.e.:
        \begin{description}
          \item
          $\bar{x}_j = 0$
          $\hspace{0.1cm} \cap \hspace{0.1cm}$
          $x^{\small LR}_j \leq \epsilon$
          $\hspace{0.2cm} \Longrightarrow \hspace{0.2cm}$
          $x_j = 0$
          \item
          $\bar{x}_j = 1$
          $\hspace{0.1cm} \cap \hspace{0.1cm}$
          $x^{\small LR}_j \geq 1 - \epsilon$
          $\hspace{0.2cm} \Longrightarrow \hspace{0.2cm}$
          $x_j = 1$
        \end{description}
      \item impose a solution time limit of $\tau$ to the optimization solver.
\end{enumerate}
A time limit is imposed since the subproblem may be difficult to solve, so the exploration of the large neighbourhood of $\bar{x},\bar{y}$ may need to be truncated. Note that in point 1 we generalize the fixing rule of RINS, in which $\epsilon = 0$. We thus allow that fixed variables may differ up to $\epsilon$ in value, in contrast to canonical RINS where they must have exactly the same value in the relaxation and in the feasible solution.

\subsection{Pheromone trail update}
At the end of each ant-construction phase $h$, we update the pheromone trails of a move $\tau_{ij}(h-1)$ according to an improved formula proposed by \cite{Ma99}:
\begin{equation}
\small
\label{pheroFormula}
\tau_{ij}(h) \hspace{0.05cm} = \hspace{0.05cm} \tau_{ij}(h-1)
\hspace{0.05cm} + \hspace{0.05cm} \sum_{k=1}^{m} \tau_{ij}^k
\hspace{0.5cm} \mbox{ with } \hspace{0.1cm} \tau_{ij}^k =
\hspace{0.1cm}
\tau_{ij}(0)
\cdot \left(
1 - \frac{z_{curr}^k - LB}{\bar{z} - LB}
\right) ,
\end{equation}

\noindent
where the values $\tau_{ij}(0)$ and $LB$ are set by using the linear relaxation of Rob-MP-CNDP: $\tau_{ij}(0)$ is set equal to the values of the corresponding optimal decision variables and $LB$ equal to the optimal value of the relaxation. Additionally, $z_{curr}^k$ is the value of the solution built by ant $k$ and $\bar{z}$ is the moving average of the values of the last $\psi$ feasible solutions built.
The formula \eqref{pheroFormula} has the desirable property of replacing the pheromone evaporation factor, a parameter whose setting may result tricky, with the moving average $\psi$, whose setting has been proved to be much less critical.

\bigskip
\noindent
Algorithm \ref{ALGhybrid} shows  the structure of our original hybrid exact-ACO algorithm. The algorithm is based on the execution of two nested loops: the outer loop is repeated until a time limit is reached and, during each execution of it, an inner loop defines $m$ ants to build the solutions. Pheromone trail updates are done at the end of each execution of the inner loop. Once the ant construction phase is over, mod-RINS is applied so to try to get an improvement by exact large neighborhood search.
\begin{algorithm}
\caption{Hybrid ACO-exact algorithm for (Rob-MP-CNDP)}
\label{ALGhybrid}
\begin{algorithmic}[1]
\State compute the linear relaxation of (Rob-MP-CNDP) and initialize the values of $\tau_{ij}(0)$ by it.
\While{the time limit is not reached}
    \For{$\mu := 1$ to $m$}
        \State build a complete routing state
        \State derive a complete feasible solution for (Rob-MP-CNDP)
    \EndFor
    \State update $\tau_{ij}(t)$ according to (\ref{pheroFormula})
\EndWhile
\State apply mod-RINS to the best feasible solution
\end{algorithmic}
\end{algorithm}

\section{Experimental results}
\label{sec:computations}

In order to assess the performance of our hybrid algorithm, we executed computational tests on a set of 30 instances, based on realistic network topologies from the \cite{SNDlib} and defined in collaboration with industrial partners from former and ongoing industrial projects (see e.g., \citealt{BlDAHa11,BlDAKa13,BaEtAl14}). The 30 instances consider 10 network topologies, whose main features are presented in Table \ref{table:instance}. For each instance, in Table \ref{table:instance} we report its identification code (\emph{ID}) and features ($|V|$ = no. vertices, $|E|$ = no. edges, $|C|$ = no. commodities). For each network topology, we defined three instances of Rob-MP-CNDP, each considering a distinct number of time periods, namely 5, 7, and 10. We performed the experiments on a machine with a 2.40 GHz quad-core processor and 16 GB of RAM and using the mixed-integer commercial solver IBM ILOG CPLEX version 12.4.
\begin{table}
\caption{Features of the instances}
\label{table:instance}
\small
\begin{center}
\begin{tabular}{c c c c c}
\hline
\hline
Name & ID & $|V|$ & $|E|$ 	& $|C|$
\\ [2pt] \hline \hline
Germany50 	& I1 & 50 & 88 & 662\\
Pioro40 	& I2 & 40 & 89 & 780\\
Norway 		& I3 & 27 & 51 & 702\\
Geant 		& I4 & 22 & 36 & 462\\
France 		& I5 & 25 & 45 & 300\\
Dfn-Gwin 	& I6 & 11 & 47 & 110\\
Pdh 		& I7 & 11 & 34 & 24\\
Ta1 		& I8 & 24 & 55 & 396\\
Polska 		& I9 & 12 & 18 & 66\\
Cost266 	& I10 & 37 & 57 & 1337\\
\hline
\hline
\end{tabular}
\end{center}
\end{table}

All the instances lead to (very) large and hard to solve Rob-MP-CNDP. We observed that even a state-of-the-art solver like CPLEX had big difficulties in identifying good feasible solutions and in the majority of cases the final optimality gap was over 95\%. In contrast, as clear from Table \ref{table:results}, in most cases our hybrid primal heuristic was able to find very high quality solutions associated with very low optimality gaps. The optimality gap indicates how far the best feasible solution found of value $v^{*}$ is from the best lower bound $LB$ available on the optimal value (formally $gap\% = |v^{*} - LB| / v^{*} \cdot 100$).
We note that using CPLEX is essential for providing a lower bound in our problem instances and thus a quality guarantee for the given solutions. In the case of solutions not produced by CPLEX, we computed the optimality gap referring to the lower bound produced by CPLEX.

On the basis of preliminary tests, we found that an effective setting of the parameters of our heuristic was: $\alpha = 0.5$ (we thus balanced attractiveness and trail level), $m = 10000$ ants, $\psi = m/10$ (width of the moving average equal to 1/10 the number of ants), $\epsilon = 0.1$ (tolerance of fixing in mod-RINS), $T = 30$ minutes (time limit imposed to the execution of mod-RINS). Each commodity admits 5 feasible paths, i.e. $|P_c| = 5, \forall c \in C$ and 5 deviations bands (2 positive, 2 negative and the null deviation band). In contrast to our first computational experience presented in \citep{DAKrPu14}, where linear relaxations were solved exactly by CPLEX requiring a non-negligible amount of time, we now 1) compute the linear relaxations of nominal problems by the results presented in Section \ref{sec:MP-CNDP} and 2) we solve the linear relaxation of Rob-MP-CNDP by a time-truncated primal simplex method as implemented by CPLEX. This allowed to greatly reduce the time of execution of the ant construction phase and hugely increase the number of defined ants.

The complete set of results is presented in Table \ref{table:results}, where we show the performance of the hybrid solution approach, that is denoted by the three measures $c^{*}$(ACO), $c^{*}$(ACO+RINS), gapAR\%, which respectively represent the value of the best solution found by pure ACO, the value of the best solution found by ACO followed by RINS and the corresponding final optimality gap. Moreover, we show the performance of CPLEX, which is denoted by measures $c^{*}$(IP) and gapIP\% representing the value of the best solution found and the corresponding final optimality gap. The value of the best solutions found for each instance is highlighted in bold type. The overall time limit for the execution of the heuristic was 1 hour. The same time limit was imposed on CPLEX when used to solve the robust counterpart Rob-MP-CNDP.  We stress that increasing the time limit did not bring any remarkable benefit to CPLEX: even when letting CPLEX run for many hours, the only effects were getting negligible improvements in the best lower bound (we observed an extremely slow improvement rate of the bound) and running out of memory because of the huge size of the search trees generated. In the case of two instances of I10, CPLEX was not even able to find a feasible solution within the time limit (cases denoted by $*$).

The best solutions found by our hybrid algorithm have in most cases a value that is at least one order of magnitude better than those found by CPLEX (3700\% better on average, excluding of course the cases for which CPLEX did not find any feasible solutions). The results are of very high quality and, given the very low optimality gap, we can suppose that some of these solutions are actually optimal. The much higher performance is particularly evident for instances I5, I8 and I9.
We observe that, in contrast to our conference paper \citep{DAKrPu14} where RINS was able to improve the value of the best solution found by the ant construction phase, in these new computations the role of RINS was reduced: thanks to the dramatic speed-up in solving the linear relaxations, we were able to really implement a swarm exploration of the feasible set, getting solutions that we believe to be optimal or very close to the optimum. So finding an improvement through RINS was impossible or really unlikely. We also think that the fact that a powerful local search like RINS has a minor or null role in finding higher quality solutions is an indication that we have defined a very strong and effective ant construction phase: in contrast to common experience in ACO, we indeed do not really require a final daemon/local search phase to get solution of sufficiently high quality.

If we focus on a specific network topology, it may be noticed that the optimality gaps produced by our algorithm tend to become smaller as the number of time periods grows. This may result surprising and counterintuitive at first sight, as one would expect exactly the opposite behaviour. Yet, considering our assumption of non-increasing demands for each commodity, it is reasonable to think that our primal hybrid heuristic discovers converging near-optimal paths when the number of time periods is large, thus performing better. We observe that it is likely that this behaviour would not occur in the case that the non-increasing demand assumption is dropped.

As further way to assess the validity and effectiveness of our new heuristic, we also used as benchmark a fast and easy method to define a possibly good feasible solution to Rob-MP-CNDP, which is suggested by the theoretical results presented in Section \ref{sec:MP-CNDP}: in each time period, we simply route the entire flow of each commodity on its shortest path. However, as clear from Table \ref{table:results}, where we report the value $c(SP)$ of the solution found by this approach and the corresponding optimality gap $gapSP\%$, in the vast majority of the instances
(in particular instances I4, I6, I7 and I8, for which the performance of our algorithm is far better),
this simple approach performed much worse than our heuristic and we believe that this provides further evidence of the solidity of our algorithm.
\begin{table}
\caption{Experimental results}
\label{table:results}
\small
\begin{center}
\begin{tabular}{c c | c c c | c c | c c}
\hline
\hline
ID 	& $|T|$ & $c^{*}$(ACO) 	& $c^{*}$(ACO+RINS) 	& gapAR\% 	 & $c^{*}(SP)$	&  gapSP\%   	 &$c^{*}$(IP) 	& gapIP\%
\\ [2pt]
\hline
\hline
& 5 & \textbf{4.26E+06} & \textbf{4.26E+06} & 3.8	&  4.35E+06    	&      5.9  	&  3.42E+08	& 98.8	 \\
I1 	&  7 	&  	\textbf{9.46E+06}	&  	\textbf{9.46E+06}	& 	2.5	&  9.64E+06  	 &      4.3  	&  8.28E+08	& 98.9	 \\
&  	  10 	&  	\textbf{3.07E+07}	&  	\textbf{3.07E+07}	& 	2.0	&  3.10E+07   	&      3.0  	&  2.89E+09	& 99.0	 \\
\hline
&  	  5 	&  	\textbf{6.14E+06}	&  	\textbf{6.14E+06}	&      11.6	&  7.48E+06    	&     27.4    	 &  2.58E+08	& 97.9
\\
I2 & 7 	& 	\textbf{1.34E+07}	& 	\textbf{1.34E+07}	&      11.5	&  1.64E+07    	&     27.6    	&  6.19E+08	& 98.1
\\
&  	10 	& 	\textbf{4.32E+07}	& 	\textbf{4.32E+07}	&       8.4	&  5.45E+07    	&     27.3    	&  2.13E+09	& 98.1	 \\
\hline
		&  	5 	& 	\textbf{3.14E+06}	& 	\textbf{3.14E+06}	&      14.6	&  3.61E+06    	 &     25.6    	 &  9.56E+07	 & 97.2
\\
I3 		&  7 	& 	\textbf{6.39E+06}	& 	\textbf{6.39E+06}	&       6.5	&  7.45E+06    	 &     19.7    	 &  2.29E+08	 & 97.4
\\
&  	10 	& 	\textbf{2.02E+07}	& 	\textbf{2.02E+07}	&       2.2	&  2.37E+07   	&     16.8    	&  7.88E+08	& 97.5	 \\
\hline
&  	5 	& 	\textbf{8.98E+05}	& 	\textbf{8.98E+05}	&       2.6	&  2.09E+06	&     58.0    	 &  4.13E+07	& 97.9	 \\
I4 		& 7 	& 	\textbf{2.15E+06}	& 	\textbf{2.15E+06}	&      11.6	&  4.94E+06  	 &     61.6    	 &  9.58E+07	 & 98.0	\\
&  	10 	& 	\textbf{6.49E+06}	& 	\textbf{6.49E+06}	&       8.8	&  1.98E+07  	&     70.1    	 &  4.18E+08	& 98.6	 \\\hline
		& 5 	& 	\textbf{1.34E+05}	& 	\textbf{1.34E+05}	&       3.0	&  1.35E+05  	 &      3.7    	 &  5.09E+06	 & 97.5	\\
I5 		& 7 	& 	\textbf{2.96E+05}	& 	\textbf{2.96E+05}	&       1.4	&  2.97E+05  	 &      1.7   	 &  1.56E+07	 & 98.1	\\
		& 10 	&  	9.43E+05	&  	9.43E+05	&       0.7	&  \textbf{9.42E+05}  	&      0.6   	 &  4.46E+07	 & 97.9	\\
\hline
		& 5 	& 	\textbf{3.34E+05}	& 	\textbf{3.34E+05}	&      12.2	&  4.14E+05  	 &     29.1    	 &  7.07E+05	 & 58.5	\\
I6 	& 7 	& 	\textbf{6.54E+05}	& 	\textbf{6.54E+05}	&       5.6	&  8.37E+05  	&     26.2    	 &  6.69E+05	&  7.7	 \\
		& 10 	& 	\textbf{2.08E+06}	& 	\textbf{2.08E+06}	&       1.7	&  2.77E+06  	 &     25.9    	 &  9.84E+07	 & 97.9	\\
\hline
		& 5 	& 	\textbf{1.19E+08}	& 	\textbf{1.19E+08}	&       1.1	&  1.50E+08  	 &     21.7    	 &  \textbf{1.19E+08}	 &  1.1	\\
I7		& 7 	& 	\textbf{2.71E+08}	& 	\textbf{2.71E+08}	&       0.5	&  3.64E+08  	 &     25.8    	 &  2.73E+08	 &  0.9	\\
		& 10 	& 	\textbf{8.29E+08}	& 	\textbf{8.29E+08}	&       0.2	&  1.27E+09  	 &     34.8    	 &  8.34E+08	 &  0.7	\\
\hline
		& 5 	& 	1.45E+08	& 	\textbf{1.39E+08}	&       9.0	&  1.60E+08  	&     21.2    	 &  7.02E+09	 & 98.2	\\
I8 		& 7 	& 	\textbf{2.83E+08}	& 	\textbf{2.83E+08}	&       3.8	&  3.36E+08  	 &     19.1    	 &  1.86E+10	 & 98.5	\\
		& 10 	& 	\textbf{9.14E+08}	& 	\textbf{9.14E+08}	&       2.1	&  1.08E+09  	 &     17.2    	 &  6.44E+10	 & 98.6	\\
\hline
		& 5 	& 	\textbf{2.12E+05}	& 	\textbf{2.12E+05}	&       1.1	&  2.22E+05  	 &      5.5   	 &  5.13E+05	 & 59.2	\\
I9 		& 7 	& 	\textbf{4.80E+05}	& 	\textbf{4.80E+05}	&       0.3	&  5.02E+05  	 &      4.9   	 &  6.08E+05	 & 21.5	\\
		& 10 	& 	\textbf{1.60E+06}	& 	\textbf{1.60E+06}	&       0.1	&  1.68E+06  	 &      4.9   	 &  9.80E+06	 & 83.7	\\
\hline
		& 5 	& 	1.37E+08	& 	1.37E+08	&      20.7	&  \textbf{1.36E+08}  	&     20.1    	 &  5.40E+09	 & 98.0	\\
I10 	& 7 	& 	\textbf{2.78E+08}	& 	\textbf{2.78E+08}	&      12.7		&  3.02E+08  &    19.8     	 &  	*	& *	 \\
		& 10 	& 	\textbf{9.26E+08}	& 	\textbf{9.26E+08}	&      14.9		&  9.84E+08  &    19.8     	 &  *		 & *	 \\
\hline
\hline
\end{tabular}
\end{center}
\end{table}

\section{Conclusion and future work}  \label{sec:end}

In this paper, we have introduced the first Robust Optimization model to handle the uncertainty that affects traffic demands in the Multiperiod Capacitated Network Design Problem (MP-CNDP). Data uncertainty may compromise the quality and feasibility of produced solutions, so, as a remedy, we have proposed a Multiband Robustness model, following the well-established methodology of Robust Optimization. We have thus produced robust solutions that are protected against deviations of input traffic data. In general, the MP-CNDP already constitutes a challenging problem, even for state-of-the-art commercial solvers like CPLEX. Accounting for robustness and considering multiple time periods has the effect of further increasing the complexity of the problem. As a matter of fact, solutions found by CPLEX are of low quality and associated with very large optimality gaps. To overcome these difficulties, we have defined a hybrid primal heuristic based on the combination of a randomized fixing algorithm inspired by ant colony optimization and an exact large neighborhood search. Computational experiments on a set of 30 realistic instances from the SNDlib have confirmed that our heuristic drastically outperforms CPLEX, finding high quality solutions associated with low optimality gaps in a short amount of time.
We believe that many of the best solutions found are actually optimal, so a future objective will be to characterize appropriate families of valid inequalities for the problem, in the attempt to close the gaps and thus possibly prove the optimality of the found solutions.
Furthermore, the excellent computational performance suggests the possibility of using the heuristic, conveniently adapted, for other applications and in more general settings, where, for example, the paths of each commodity are not predetermined.
We expect our hybrid primal heuristic to perform well even in such different contexts.

\section*{Acknowledgements}
This work was partially supported by the \emph{German Research Foundation} (DFG), project \emph{Multiperiod Network Optimization}, by the DFG Research Center MATHEON (www.matheon.de), Project B3, by the Einstein Campus for Mathematics Berlin (ECMath) Project ROUAN, and by the \emph{German Federal Ministry of Education and Research} (BMBF), Project \emph{ROBUKOM} \citep{BaEtAl14}, grant 05M10PAA.

\appendix

\section{Proofs of the Statements of Section \ref{sec:MP-CNDP}}
\label{Appendix}

\subsection{Proof of Proposition \ref{Pr:consecutivePeriodFlow}}

\noindent
\emph{
\textbf{Statement:}
Let $\bar{x}^t_{c p}:=0$, $d^t_{c p}:=0$ and $\gamma^t_{e}:=\infty$ for $t=0$. There exists an optimal solution $\bar{x}$ of MP-CNDP-LP such that:
\begin{equation*}
 d^{t-1}_c \hspace{0.1cm} \bar{x}^{t-1}_{c p}
\hspace{0.15cm} \leq \hspace{0.15cm}
d^t_c \hspace{0.1cm}  \bar{x}^t_{c p}
\hspace{1.0cm}
\forall \hspace{0.1cm} c \in C, p \in P_c, t \in T
\end{equation*}
}

\begin{Proof}
We prove the statement by induction, considering two consecutive time periods $t-1$ and $t$ with $t = \{1,2,\ldots,|T|\}$. As basis step, we note that relation \eqref{eq:trafficRelation} trivially holds for $t = 1$ by definition of $\bar{x}^t_{c p}:=0$, $d^t_{c p}:=0$ for $t = 0$.
As inductive step, suppose that relation \eqref{eq:trafficRelation} holds for $t \leq |T| - 1$. We will show that it holds also for $t \leq |T|$ period. Suppose that relation \eqref{eq:trafficRelation} does not hold for $t = |T|$, i.e. $\exists \hspace{0.1cm} c^{*} \in C, p^{*} \in P_{c^{*}}:$ $d^{t-1}_{c^{*}} \bar{x}^{t-1}_{c^{*} p^{*}} > d^t_{c^{*}} \bar{x}^{t}_{c^{*} p^{*}}$ for some commodity and its corresponding path when $t=T$. Without loss of generality, we assume that there exists unique $c^{*} \in C, p^{*} \in P_{c^{*}}$ that does not respect relation \eqref{eq:trafficRelation} for period $|T|$, as the argument we use below can be applied iteratively to any number of such deviating inequalities. Moreover, suppose that $p^{*} \in P_c$
is a shortest path in edge cost for $c^{*}$. Let $\Delta = d^{t-1}_{c^{*}}  \bar{x}^{t-1}_{c^{*} p^{*}} - d^t_{c^{*}} \bar{x}^{t}_{c^{*} p^{*}}$ be the residual flow.  For simplicity we assume that $\Delta$ is routed over  a single path as the argument similarly holds when it splits over several paths. Given our assumption, it is easy to see that $\Delta$ will be routed on another path $p' \in P_{c^{*}}$ and that the objective function value will gain at most $\sum\limits_{e\in E: e \in p'}\gamma^t_e\Delta \geq \sum\limits_{e\in E: e \in p^{*}}\gamma^t_e\Delta$. When $p'$ is another shortest path, we get an equality and we see that $x^{t-1}_{c^{*} p^{*}} d^{t-1}_{c^{*}} \leq x^t_{c p^{*}} d^t_{c^{*}}$ gives an equivalent optimal solution in terms of the objective function value, otherwise we have a contradiction.  Suppose $p^{*}$ in  $x^{t}_{c^{*} p^{*}}$  is not  a shortest path in edge cost for $c^{*}$. If  $\sum\limits_{e\in E: e \in p^{*}}\gamma^t_e\leq \sum\limits_{e\in E: e \in p'}\gamma^t_e$ then the same argument as above applies. If  $\sum\limits_{e\in E: e \in p^{*}}\gamma^t_e > \sum\limits_{e\in E: e \in p'}\gamma^t_e$ then the solution gains at most $\sum\limits_{e\in E: e \in p^{*}}\gamma^t_e\Delta$ $-$  $\sum\limits_{e\in E: e \in p'}\gamma^t_e\Delta > 0$. If this is the case, we can reroute $\Delta$ in period $|T|-1$ from $p^{*}$ to $p'$. Then we have $\bar{x}^{t-1}_{c^{*} p^{*}}d^{t-1}_{c^{*}} = \bar{x}^t_{c^{*} p^{*}}d^t_{c^{*}}$ in period T, when in the deviating inequality $\bar{x}^{t-2}_{c^{*} p^{*}}d^{t-2}_{c^{*}} \leq \bar{x}^{t}_{c^{*} p^{*}}d^{t}_{c^{*}}$ holds. Then this new solution gains at most $\sum\limits_{e\in E: e \in p^{*}}\gamma^{t-1}_e\Delta$ $-$  $\sum\limits_{e\in E: e \in p'}\gamma^{t-1}_e\Delta \geq \sum\limits_{e\in E: e \in p^{*}}\gamma^t_e\Delta$ $-$  $\sum\limits_{e\in E: e \in p'}\gamma^t_e\Delta$ by the non-increasing cost per unit of capacity. The same argument holds when $\bar{x}^{t-2}_{c^{*} p^{*}}d^{t-2}_{c^{*}} > \bar{x}^{t}_{c^{*} p^{*}}d^{t}_{c^{*}}$.
\qed
\end{Proof}

\bigskip

\subsection{Proof of Theorem \ref{Th:optSolMP-CNDP-LP}}

\noindent
\emph{
\textbf{Statement:}
Consider problem MP-CNDP-LP, namely the linear relaxation of MP-CNDP-IP, and let $p_c^{*} \in P_c$ be a shortest path in edge length for each commodity $c \in C$. An optimal solution of MP-CNDP-LP can be defined by routing in each time period $t \in T$ the entire flow of each commodity $c \in C$ on the shortest path $p_c^{*}$ and by installing on each edge of $p_c^{*}$ the exact capacity needed to route the traffic flow $d_{ct}$.
}

\begin{Proof}
Before proceeding to the proof, we recall that in MP-CNDP-LP the capacity variables lose their integrality constraint and thus the capacity installation on each edge in each period will exactly equal the flow.
It is well known that the result stated in the theorem holds for one single period (see \citealt{AhMaOr93}) and we propose a proof of the result, extending below the reasoning to the case with multiple periods.
Consider then a single-period problem (in this case, the $t$ apex is thus equal to 1 for all the involved quantities) and a commodity $c \in C$. Suppose that $x^1_{c p}=k$ with $0 < k \leq 1$ is a non-zero variable in the optimal solution of the problem for some $p \in P_c$. Hence the contribution to the objective function is $\sum\limits_{e\in E: e \in p}\gamma^{1}_e d^1_{c} \hspace{0.1cm} k$. Let $p'$  be a shortest path for $c$. Then $\sum\limits_{e\in E: e \in p'}\gamma^1_e  d^1_{c} \hspace{0.1cm}( x^1_{c p'}+k) \leq \sum\limits_{e\in E: e \in p'}\gamma^1_e  d^1_{c} \hspace{0.1cm} x^1_{c p'} + \sum\limits_{e\in E: e \in p}\gamma^1_e d^1_{c} \hspace{0.1cm} x^1_{cp}$ . Thus the theorem holds for the single period case.
Assume now that the theorem holds for $|T|-1$ periods. We will show that it holds for $|T|$ periods. By Proposition \ref{Pr:consecutivePeriodFlow}, we may assume that  $x^{t-1}_{c p}d^{t-1}_c \leq x^t_{c p}d^t_c$ for all $c \in C, p \in P_c, t \in T$ in an optimal solution of MP-CNDP-LP with $|T|$ periods.  Let $\overline{OPT}(|T|-1)$
denote the optimal objective function value for $|T|-1$ periods according to the induction hypothesis. Also, let $OPT_t$ denote the contribution  of the optimal solution of the  MP-CNDP-LP for $|T|$ periods to the objective function in period $t$.
Then $\sum\limits_{t \in T} OPT_t = \sum\limits_{t \in T \setminus |T|} OPT_t  + OPT_{|T|}$. Suppose $\sum\limits_{t \in T \setminus |T|} OPT_t$ is determined by, for each commodity, the same shortest path in edge cost for each period up to $|T|-1$. Hence $\overline{OPT}(|T|-1) = \sum \limits_{t \in T \setminus |T|} OPT_t $. Using a similar argument as in the single period case, it is easy to see that the increase in demand for each commodity in period $|T|$ will be routed in the shortest paths that were used in the previous periods. For this case, let $\overline{OPT}_{|T|}$  denote the contribution of the optimal solution of MP-CNDP-LP to the objective function in period $|T|$.  Suppose $\sum \limits_{t \in T \setminus |T|} OPT_t$ is not determined by, for each commodity, the same shortest path in edge cost for each period up to $|T|-1$. Then it is clear that $\overline{OPT}(|T|-1) \leq \sum\limits_{t \in T \setminus |T|} OPT_t $ and $\overline{OPT}(|T|)\leq OPT_{|T|}$. Hence $\overline{OPT}(|T|-1) + \overline{OPT}(|T|) \leq  \sum\limits_{t \in T \setminus |T|} OPT_t +  OPT_{|T|}$.
\qed
\end{Proof}


\section*{References}


\begin{thebibliography}{00}

\bibitem[Ahuja et al.(1993)]{AhMaOr93}
Ahuja, R.K., Magnanti, T., Orlin, J.: Network Flows: Theory, Algorithms, and Applications. Prentice Hall, Upper Saddle River (1993)

\bibitem[Bauschert et al.(2014)]{BaEtAl14}
Bauschert, T., B\"using, C., D'Andreagiovanni, F., Koster, A.M.C.A., Kutschka, M., Steglich, U.: Network Planning under Demand Uncertainty with Robust Optimization. IEEE Communications Magazine, \textbf{52} (2), 178-185 (2014)
DOI: 10.1109/MCOM.2014.6736760

\bibitem[Belotti et al.(2011)]{BeKoNo11}
Belotti, P., Kompella, K., Noronha, L.: A comparison of OTN and MPLS networks
under traffific uncertainty. Submitted in 2011.
http://myweb.clemson.edu/$\sim$pbelott/papers/robust-opt-network-design.pdf
(retrieved on 07.03.2014) (2011)

\bibitem[Ben-Tal et al.(2009)]{BeElNe09}
Ben-Tal, A., El Ghaoui, L., Nemirovski, A.: Robust Optimization. Springer, Heidelberg (2009)

\bibitem[Bertsekas(1998)]{Be98}
Bertsekas, D.P.: Network Optimization: Continuous and Discrete Models. Athena Scientific, Belmont (1998)

\bibitem[Bertsimas et al.(2011)]{BeBrCa11}
Bertsimas, D., Brown, D., Caramanis, C.: Theory and Applications of Robust Optimization. SIAM Review \textbf{53} (3), 464-501 (2011)

\bibitem[Bertsimas and Sim(2004)]{BeSi04}
Bertsimas, D., Sim, M.: The Price of Robustness. Oper. Res. \textbf{52} (1), 35-53 (2004)

\bibitem[Bienstock et al.(1998)]{BiEtAl98}
Bienstock, D., Chopra, S., G\"unl\"uk, O, Tsai, C.: Minimum cost capacity installation for multicommodity network flows. Math. Prog. \textbf{81} (2), 177-199 (1998)

\bibitem[Bley et al.(2011)]{BlDAHa11}
Bley, A., D'Andreagiovanni, F, Hanemann, D.: Robustness in Communication Networks: Scenarios and Mathematical Approaches. In: Proc. of the ITG Symposium on Photonic Networks 2011, 1-8, VDE Verlag, Berlin (2011)

\bibitem[Bley et al.(2013)]{BlDAKa13}
Bley, A., D'Andreagiovanni, F, Karch, D.: WDM Fiber Replacement Scheduling. In: Proc. of INOC 2013,  Electronic Notes in Discrete Mathematics \textbf{41} (5), 189-196 (2013)
DOI: 10.1016/j.endm.2013.05.092

\bibitem[Bley et al.(2000)]{BlGrWe07}
Bley, A., Gr\"otschel, M., Wess\"aly, R.: Design of broadband virtual private networks: Model and heuristic for the B-WiN. DIMACS Series in Discrete Mathematics and Theoretical Computer Science
53, 1-16 (2000)

\bibitem[Blum (2005)]{Bl05}
Blum, C.: Ant colony optimization: Introduction and recent trends. Physics Life Rew. \textbf{2} (4), 353-373 (2005)

\bibitem[Blum et al.(2011)]{BlEtAl11} Blum., C., Puchinger, J., Raidl, G.R., Roli., A.: Hybrid metaheuristics in combinatorial optimization: A survey. Appl. Soft Comp. \textbf{11} (6), 4135--4151 (2011)

\bibitem[B\"using and D'Andreagiovanni(2012)]{BuDA12a} B\"using, C., D'Andreagiovanni, F.: New Results about Multi-band Uncertainty in Robust Optimization. In: Klasing, R. (ed.) Experimental Algorithms - SEA 2012, LNCS, vol. 7276, pp. 63-74. Springer, Heidelberg (2012)
    DOI: 10.1007/978-3-642-30850-5\_7

\bibitem[B\"using and D'Andreagiovanni(2013)]{BuDA13} B\"using, C., D'Andreagiovanni, F.: Robust Optimization under Multi-band Uncertainty - Part I: Theory. CoRR abs/1301.2734, http://arxiv.org/abs/1301.2734 (2013)

\bibitem[B\"using and D'Andreagiovanni(2014)]{BuDA14} B\"using, C., D'Andreagiovanni, F.: A new theoretical framework for Robust Optimization under multi-band uncertainty. In: Helber, S. et al. (eds.) Operations Research Proceedings 2012, pp. 115-121. Springer, Heidelberg (2014)

\bibitem[B\"using et al.(2014)]{BuDARa14}
B\"using, C., D'Andreagiovanni, F., Raymond, A.: 0-1 Multiband Robust Optimization. In: Huisman, D. et al. (eds.) Operations Research Proceedings 2013, pp. 89-95. Springer, Heidelberg (2014)

\bibitem[Christofides and Brooker(1974)]{ChBr74}
Christofides, N., Brooker, P.: Optimal Expansion of an Existing Network. Math. Program. \textbf{6}, 197-211 (1974)

\bibitem[Crainic and Gendreau(2002)]{CrGe02}
Crainic, T.G., Gendreau, M.: Cooperative Parallel Tabu Search for Capacitated Network Design. J. Heuristics \textbf{8} (6), 601-627 (2002)

\bibitem[Danna et al.(2005)]{DaRoLP05}
Danna, E., Rothberg, E., Le Pape, C.: Exploring relaxation induced neighborhoods to improve MIP solutions. Math. Program. \textbf{102}, 71-90 (2005)

\bibitem[D'Andreagiovanni et al.(2014)]{DAKrPu14}
D'Andreagiovanni, F., Krolikowski, J., Pulaj, J.: A hybrid primal heuristic for Robust Multiperiod Network Design. To appear in Esparcia-Alcazar, A. (ed.) Applications of Evolutionary Computation - 17th European Conference, EvoApplications 2014, Lecture Notes in Computer Sciences. Springer, Heidelberg (2014)

\bibitem[Dely et al.(2013)]{DeDAKa13} Dely, P., D'Andreagiovanni, F., Kassler, A.: Fair optimization of mesh-connected WLAN hotspots. Wireless Communications and Mobile Computing (2013) DOI: 10.1002/wcm.2393

\bibitem[Dorigo et al.(1999)]{DoDCGa99}
Dorigo, M., Di Caro, G., Gambardella, L.M.: Ant algorithms for discrete optimization.
Artificial Life \textbf{5} (2), 137-172 (1999)

\bibitem[Dorigo et al.(1996)]{DoMaCo96}
Dorigo, M., Maniezzo, V., Colorni, A.: Ant System: Optimization by a colony of cooperating agents.
IEEE Trans. Syst. Man Cybern. B \textbf{26} (1), 29-41 (1996)

\bibitem[Gambardella et al.(2012)]{GaMoWe12}
Gambardella, L.M., Montemanni, R., Weyland, D.: Coupling ant colony systems with strong local searches. Europ. J. Oper. Res. \textbf{220} (3), 831-843 (2012)

\bibitem[Gavros and Raghavan(2012)]{GaRa12}
Gamvros, I., Raghavan, S.: Multi-period traffic routing in satellite networks. Europ. J. Oper. Res. \textbf{219} (3), 738-750 (2012)

\bibitem[Gendrau et al.(2006)]{GeEtAl06}
Gendreau, M., Potvin, J., Smires A., Soriano P.: Multi-period capacity expansion for a local access telecommunications network. Europ. J. Oper. Res., \textbf{172} (3), 1051-1066 (2006)

\bibitem[Gupta and Grossmann(2012)]{GuGr12}
Gupta, V., Grossmann, I.E.: An Efficient Multiperiod MINLP Model for Optimal Planning of Offshore Oil and Gas Field Infrastructure. Ind. Eng. Chem. Res., \textbf{51} (19), 6823–6840 (2012)

\bibitem[CPLEX(2014)]{CPLEX} IBM ILOG CPLEX,
        {\small http://www-01.ibm.com/software/integration/optimization/cplex-optimizer.}

\bibitem[Kleeman et al.(2012)]{KlEtAl12}
Kleeman, M.P., Seibert, B.A., Lamont, G.B., Hopkinson, K.M., Graham, S.R.: Solving Multicommodity Capacitated Network Design Problems Using Multiobjective Evolutionary Algorithms. IEEE Trans. Evol. Comp. \textbf{16} (4), 449-471 (2012)

\bibitem[Koster et al.(2013)]{KoKuRa13}
Koster, A.M.C.A., Kutschka, M., Raack, C.: Robust network design: Formulations, valid inequalities, and computations. Networks \textbf{61} (2), 128-149 (2013)

\bibitem[Lardeux et al.(2007)]{LaNaGe07}
Lardeux, B., Nace, D., Geffard, J.: Multiperiod network design with incremental routing.
Networks \textbf{50} (1), 109-117 (2007)

\bibitem[Maniezzo(1999)]{Ma99}
Maniezzo, V.: Exact and Approximate Nondeterministic Tree-Search Procedures for the Quadratic Assignment Problem. {\em INFORMS J. Comp.} \textbf{11} (4), 358-369 (1999)

\bibitem[SNDlib(2010)]{SNDlib}
Orlowski, S., Wess\"aly, R., Pioro, M., Tomaszewski, A.: SNDlib 1.0 - Survivable Network Design Library. Networks \textbf{55} (3), 276-286 (2010)

\bibitem[Theimer(2010)]{Th09}
Theimer, T: Towards a scalable and flat IP core network. {\em ECOC 2009} Vienna, Austria.
(retrieved 07.03.2014)
\\
http://conference.vde.com/ecoc-2009/programs/documents/ws4\_thomas\%20theimer.pdf

\end{thebibliography}
\end{document}